\documentclass[11pt]{article}
\usepackage{amsmath, amsbsy, amscd, amsthm, amsfonts, array, epic, eepic, verbatim}
 \newtheorem{thm}{Theorem}
 
 \newtheorem{lem}[thm]{Lemma}
 
 \theoremstyle{definition}

 \numberwithin{equation}{section}

 \newcommand{\B}{\mathcal{B}}

 \newcommand{\OO}{\mathcal{O}}

 \newcommand{\X}{\mathcal{X}}

 \newcommand{\LHS}{{\rm LHS}}
 \newcommand{\RHS}{{\rm RHS}}

\def\ZZ{\mathbb{Z}}
\def\FF{\mathbb{F}}
\def\QQ{\mathbb{Q}}

\def\CC{\mathbb{C}}

\def\EE{\mathbb{E}}
\def\AA{\mathbb{A}}
\def\PP{\mathbb{P}}
\def\ov{\overline}

\def\Ker{\operatorname{Ker}}

\def\Sym{\operatorname{Sym}}

\def\Hom{\operatorname{Hom}}

\def\det{\operatorname{det}}

\begin{document}  

\title{ \Large{THE CREPANT RESOLUTION CONJECTURE FOR\\ 3-DIMENSIONAL FLAGS MODULO AN INVOLUTION} }
\author{W.~D.~Gillam \\ Columbia University Department of Mathematics}

\date{June 6, 2007}

\maketitle

\begin{abstract} After fixing a non-degenerate bilinear form on a vector space $V$ we define a $\ZZ_2$-action on the manifold of flags $F$ in $V$ by taking a flag to its orthogonal complement.  When $V$ is of dimension 3 we check that the Crepant Resolution Conjecture of J.~Bryan and T.~Graber holds: the genus zero (orbifold) Gromov-Witten potential function of $[F / \ZZ_2]$ agrees (up to unstable terms) with the genus zero Gromov-Witten potential function of a crepant resolution $Y$ of the quotient scheme $F / \ZZ_2$, after setting a quantum parameter to $-1$, making a linear change of variables, and analytically continuing coefficients.  We explicitly compute all degree $0$ and $3$-point invariants for the orbifold and the resolution, then argue that the other invariants are determined by WDVV and the Divisor Axiom.  The quotient $F / \ZZ_2$ is contained in the quotient of $\PP^2 \times \PP^2$ by the $\ZZ_2$-action interchanging the factors and the crepant resolution $Y$ (a hypersurface in the Hilbert scheme ${\rm Hilb}^2 \PP^2$) is the projectivization of a novel rank 2 vector bundle over $\PP^2$.
\end{abstract}

\begin{section}{Introduction} \label{section:intro}

Fix a non-degenerate bilinear form $\langle , \rangle$ on an $n$-dimensional complex vector space $V$.  For a linear subspace $A \subseteq V$ let $$A^\perp := \{ v \in V : \langle v,a \rangle = 0 {\rm \; for \; all \;} a \in A \}$$ denote the ``orthogonal complement'' of $A$ with respect to $\langle, \rangle$.  This is of dimension $n -\dim A$ but may not be disjoint from $A$.  We may identify $\Hom(V/A,\CC)$ with $A^\perp$ via this inner product.  Then $$(A_1 \subset A_2 \subset \dots \subset A_{n-1}) \mapsto ( A_{n-1}^\perp \subset \dots \subset A_{2}^\perp \subset A_1^\perp)$$ is an involution of ($\ZZ_2$-action on) $F$, the manifold of complete flags in $V$.  For example, consider the bilinear form $$\langle v,w \rangle := v_nw_1 + v_{n-1}w_2 + \dots + v_1 w_n$$ on $\CC^n$.  Let $W$ be the idempotent $n \times n$ matrix whose $(i,j)$ entry is $\delta_i^{n+1-j}$, so that multiplying by $W$ on the right reverses the columns of a matrix and multiplying on the left reverses the rows.  The idempotent outer automorphism $A \mapsto W(A^T)^{-1}W$ of $G=SL_n(\CC)$ preserves the Borel subgroup $B^+$ of upper triangular matrices, hence induces a $\ZZ_2$-action on $F=G/B^+$ which takes a flag to its orthogonal complement with respect to the above bilinear form.

The goal of this paper is to check the Crepant Resolution Conjecture (CRC) \cite{BG} of J.~Bryan and T.~Graber when $V$ is of dimension $3$.  This conjecture asserts that the Gromov-Witten potential function of $Y$, where $Y$ is any resolution of singularities $r:Y \mapsto X / G$ satisfying $\omega_Y = r^* \omega_{ X / G }$ (i.e. a \emph{crepant} resolution), should be equivalent to the Gromov-Witten potential function of the stack (or orbifold) $\X := [X / G]$, assuming that $\X$ is a Gorenstein orbifold satisfying the Hard Lefschetz Condition (as is the case for the example considered here: $X=F$, $G=\ZZ_2$).  This orbifold potential function involves integrals over moduli spaces of stable maps from ``orbicurves" to $\X$ (see \cite{AGV}, \cite{AGV2}, \cite{CR2} for discussion).  The classes to be integrated are to be pulled back (via evaluation maps) from the Chen-Ruan orbifold cohomology ring \cite{CR} of $\X$.  As a vector space, this is the cohomology of the ``inertia stack" of $\X$ which, in our case, is simply the direct sum of the $\ZZ_2$-invariant cohomology of $F$ and the cohomology of the fixed locus $F^{\ZZ_2}$.  This vector space carries the structure of a graded ring, whose product can be deformed using $3$-point, genus zero orbifold Gromov-Witten invariants to form a quantum orbifold cohomology ring (associativity is proved in \cite{AGV2}).  The potential functions are to be identified after performing a linear (and grading preserving) change of variables, throwing away ``unstable terms" (passing, say, to third derivatives of the potential functions), and setting extra quantum parameters to roots of unity in the potential of the resolution.  The extra quantum parameters in the potential of a crepant resolution will be degree zero quantum parameters corresponding to curves on the exceptional locus of the resolution map.  One does not expect the cohomology of the resolution and the usual orbifold cohomology to agree, because these degree zero quantum parameters must first be set to roots of unity to deform the multiplication.  However, after doing this, the two quantum cohomology rings should be isomorphic; we show that this is the case for $\X$, $Y$ in Section~\ref{section:changevars}.  We use this isomorphism to work out the general change of variables and check that the genus zero potential functions agree in Section~\ref{section:higherpoint}.

In Section~\ref{section:orbcoho} we recall some general facts about flag manifolds, then discuss the $\ZZ_2$-action on the manifold of flags in $\CC^3$, determining the fixed locus and the corresponding restriction map on cohomology, leading to a description of the Chen-Ruan orbifold cohomology ring of $\X$.  Then we find explicit descriptions of some simple moduli spaces of orbifold stable maps to $\X$, which, together with associativity of the orbifold quantum product, we use to determine the 3-point Gromov-Witten invariants and the orbifold quantum cohomology ring of $\X$.  In Section~\ref{section:y} we describe the crepant resolution $r:Y \to F / \ZZ_2$ as a $\PP^1$-bundle over $\PP^2$ which is a hypersurface in the Hilbert scheme ${\rm Hilb}^2 \PP^2$.  In \cite{W}, J.~Wise showed (using Graber's computations \cite{G}) that the CRC holds for $\X=[\PP^2 \times \PP^2 / \ZZ_2]$, $Y = {\rm Hilb}^2 \PP^2$.  It may be possible to use this result together with the Orbifold Quantum Lefschetz Hyperplane Theorem (\cite{HHT}, Section~5) to compute the orbifold potential of $[F / \ZZ_2]$.  In Section~\ref{section:gwy}, we compute the 3-point Gromov-Witten invariants of $Y$, making heavy use of WDVV and we give a multiplication table and presentation of the quantum cohomology of $Y$.  

I happily credit my advisor, Michael Thaddeus, with many ideas presented in this paper. 

\end{section}

\begin{section}{Orbifold Quantum Cohomology of $\X$} \label{section:orbcoho}

Abusing notation, we let $A_i$ denote the rank $i$ vector bundle on $F$ whose fiber over $(A_1 \subset \dots \subset A_{n-1})$ is $A_i$.  Let $u_i:=c_1(A_{i+1}/A_i)$ and let $p_i := -c_1( \det A_i)$.  Applying adjunction to the short exact sequence (SES) $$0 \to A_i \to A_{i+1} \to A_{i+1}/A_i \to 0$$ we get $-p_{i+1} = u_i-p_i$ (one should put a ``dual'' or a minus sign somewhere in the fifth sentence in the second paragraph of page 3 in \cite{GK}).  Once a basis $\{ e_1, \dots, e_n \}$ for $V$ is chosen, the effective cone $H^{\rm eff}_2(F)$ of $F$ is spanned by the curve classes $$W_i := \{ (A_1 \subset \dots \subset A_{n-1}) \in F: A_j = \langle e_1,\dots,e_j \rangle {\rm \, unless \,} j=i \},$$ which satisfy $\langle W_i, p_j \rangle = \delta_{ij}$ so the $p_i$ are dual to the effective cone under the evaluation pairing.  The $\ZZ_2$-action takes $W_i$ to $W_{n-i}$ (hence $p_i$ to $p_{n-i}$ in cohomology).  Dualizing the SES $$0 \to A_i \to V \to V / A_i \to 0$$ we find that $c_1(A_i)=c_1(A_i^\perp)$.

Specialize to the case where $V=\CC^3$.  In terms of the $p_i$, the cohomology ring of $F$ can be presented: $$H^*(F,\ZZ)=\ZZ[p_1,p_2]/\langle p_1p_2-p_1^2-p_2^2,p_1p_2^2-p_1^2p_2 \rangle$$ (the relations are given by the elementary symmetric polynomials in the $u_i$).  The forgetful maps $$\begin{array}{l} \pi_1(A_1 \subset A_2) := A_1 \in \PP^2 \\ \pi_2(A_1 \subset A_2) \mapsto A_2 \in \PP^{2*} := {\rm Gr}_2 (\CC^3) \end{array}$$ make $F$ a $\PP^1$-bundle over $\PP^2$, $\PP^{2*}$ respectively.  The classes $p_1$, $p_2$ are pullbacks of the positive generators of $H^2(\PP^2,\ZZ)$, $H^2(\PP^{2*},\ZZ)$ via these projections.

The fixed locus $C := F^{\ZZ_2}$ is given by $$\begin{array}{ll} \{(A_1 \subset A_2):A_2^\perp=A_1 \} & =  \{ (A_1 \subset A_2): A_2 = A_1^\perp \} \\
 & = \{ (A_1 \subset A_1^\perp) : \langle A_1,A_1 \rangle=0 \} \\ & = \{ (A_2^\perp \subset A_2) : \langle , \rangle|_{A_2} {\rm \; is \; degenerate} \}. \end{array}$$ This is a section of both of the above projections over a conic in $\PP^2, \PP^{2*}$ so $C \cong \PP^1$ and the restriction map $H^*(F,\ZZ) \to H^*(C,\ZZ)$ takes $p_1$ and $p_2$ to twice the positive generator $x$ of $H^2(C,\ZZ) \cong \ZZ$.

Next we compute the Poincar\'e dual of $[C]$ in $H^4(F,\ZZ)$.  This is characterized by the property that $$\langle a \cup {\rm P.D.} [C] , [F] \rangle = \langle a|_C , [C] \rangle {\rm \;\; for \; all \;} a \in H^2(F,\ZZ)$$ so, since we know $\langle p_1^2p_2, [F] \rangle = \langle p_1p_2^2,[F] \rangle = 1$, we find that $${\rm P.D.} [C] = 2p_1p_2.$$  

\bigskip

\centerline {\bf Orbifold cohomology of $\X$}

\bigskip

From this we can compute the Chen-Ruan orbifold cohomology ring $H^*_{orb}( \X )$.  Notice that $p_1 + p_2$ generates $H^*(F,\QQ)^{\ZZ_2}$ as a $\QQ$-algebra.  Additively (with complex grading) we will use the basis $S_0, \dots, S_5$ below: $$\begin{array}{ccc} {\rm grading/sector} & H^*(F,\QQ)^{\ZZ_2} & H^*(C,\QQ) \\ 0 & S_0=1 \\ 1 & S_1=p_1+p_2 & S_2=1 \\ 2 & S_3 = (p_1+p_2)^2 = 3p_1p_2 & S_4=x \\ 3 & S_5=(p_1+p_2)^3=6p_1^2p_2 \end{array}$$ with multiplication $S_1S_2=(p_1+p_2)|_C=4S_4$, $S_2^2 = P.D.[C]= (2/3)S_3$, $S_2S_4=(1/6)S_5$.  Thus the Poincar\'e duality metric $G=(G_{ij})=(\langle S_iS_j, [F / \ZZ_2] \rangle)$ and the corresponding dual basis are $$G=\begin{pmatrix} 0 & 0 & 0 & 0 & 0 & 3 \\ 0 & 0 & 0 & 3 & 0 & 0 \\ 0 & 0 & 0 & 0 & 1/2 & 0 \\ 0 & 3 & 0 & 0 & 0 & 0 \\ 0 & 0 & 1/2 & 0 & 0 & 0 \\ 3 & 0 & 0 & 0 & 0 & 0 \end{pmatrix} \quad \quad \begin{array}{ll} S^0 = (1/3)S_5 & \\ S^1 = (1/3)S_3 & S^2 = 2S_4 \\ S^3 = (1/3)S_1 & S^4 = 2S_2 \\ S^5 = (1/3)S_0 \end{array}$$ and a nice presentation is $H^*_{orb}( \X ) \cong \QQ[S_1,S_2] / \langle S_2^3,3S_2^2-2S_1^2 \rangle.$

\bigskip

\centerline {\bf Orbifold stable maps to $ \X $}

\bigskip

\noindent Here we summarize the treatment of twisted stable maps given in $\S 3$ of \cite{AGV}, restricting to the case of a $\ZZ_2$ global quotient orbifold\footnote{According to some definitions, the $\ZZ_2$-action should be generically free for the stack quotient to be an orbifold, but we certainly want to allow $X$ to be a point here.} $\X=[X / \ZZ_2]$.  Assume the $\ZZ_2$-fixed locus is connected for simplicity\footnote{Otherwise the moduli spaces discussed here should be more carefully divided into components depending on which of the fixed components a given marked ramification point is mapped to.}.  For a fixed $\ZZ_2$-invariant effective homology class $\beta \in H_2^{\rm eff}(X)^{\ZZ_2}$ and integers $r,u$, let $\ov{M}_{0;r,u}(\X,\beta)$ denote the DM-stack representing (isomorphism classes of) flat families of commutative diagrams $$\begin{CD} \tilde{\Sigma} @>{\phi}>{\ZZ_2{\rm -equivariant}}> X \\ @V{q}VV @VVV \\ (\Sigma,R_1,\dots,R_r,U_1,\dots,U_u) @>>> X/ \ZZ_2 \end{CD}$$ where $\Sigma$ is a connected nodal curve (with marked points $R_i,U_i$ in its smooth locus) of arithmetic genus zero, $\tilde{\Sigma}$ is a curve (possibly disconnected of positive genus) with generically-free $\ZZ_2$-action branched over each $R_i$ and possibly over nodes of $\Sigma$, but nowhere else.  The quotient map $q$ must take nodes to nodes, $\phi$ is required to be $\ZZ_2$-equivariant, and we require $\phi_*[\tilde{\Sigma}]=\beta$.  For stability reasons, we require that any component of $\Sigma$ mapped to a point have at least $3$ special points (marked points and nodes).  For monodromy reasons, the moduli space $\ov{M}_{0;r,u}( \X , \beta)$ is empty unless $r$ is even.  The evaluation maps $\ov{e}_i$ at ramification points naturally take values in the fixed locus $X^{\ZZ_2}$ and evaluation $e_i$ at non-ramification points naturally takes values in the quotient $X/\ZZ_2$, so we can define orbifold Gromov-Witten invariants as follows.  For $\alpha_1,\dots,\alpha_r \in H^*(X^{\ZZ_2},\QQ)$ and $\phi_1, \dots , \phi_u \in H^*(X,\QQ)^{\ZZ_2}=H^*(X/ \ZZ_2,\QQ)$ we define $$\langle \alpha_1 , \dots , \alpha_r , \phi_1 , \dots , \phi_u \rangle^\beta := \int_{ \ov{M}_{0;r,u}( \X ,\beta) } (\ov{e}_1^* \alpha_1) \cdots (\ov{e}_r^* \alpha_r) (e_1^* \phi_1) \cdots (e_u^* \phi_u).$$

Returning to the case at hand $\X = [F / \ZZ_2]$, we will identify an effective $\ZZ_2$-invariant homology class $aW_1+aW_2$ with the non-negative integer $a$.  To compute the orbifold quantum cohomology of $\X$ it will suffice to study the moduli spaces $\ov{M}_{0;0,3}(\X,1)$ and $\ov{M}_{0;2,0}( \X ,1)$, which have expected dimension (this will be equal to the actual dimension) $5$ and $2$ respectively.  The quantum parameter $q$ for $QH^*_{orb}([F / \ZZ_2])$ should have degree $2$ because $$\langle c_1(T \X), 1 \rangle = (1/2) \langle c_1(F),  W_1 + W_2 \rangle = (1/2) \langle 2p_1+2p_2, W_1+W_2 \rangle = 2.$$  We may identify the first of these moduli spaces (at least coarsely) with the usual Kontsevich stable map space $\ov{M}_{0;3}(F,(1,0))$ by the map $$\begin{array}{c}[ f:(\Sigma,P_1,P_2,P_3) \to F ] \mapsto [f \coprod f^\perp : \tilde{\Sigma} := \Sigma \coprod \Sigma \to F] \end{array}$$ in which case the evaluation maps $e_i : \ov{M}_{0;0,3}( \X ,1) \to F / \ZZ_2$ are identified with $ev_i+ev_i^\perp$.  Thus we have a commutative diagram $$\begin{CD} \ov{M_{0;3}}(F,(1,0)) @>{ev_i}>> F \\ @V{\cong}VV @V{}VV \\ \ov{M}_{0;0,3}( \X ,1) @>{e_i}>> F / \ZZ_2 \end{CD}$$ from which we can compute $\langle S_1,S_3,S_3 \rangle^1  =  9$ and $\langle S_1,S_1,S_5 \rangle^1  = 6$ (c.f. \cite{GK}).

Next we claim that the evaluation map $$\ov{e} = (\ov{e}_1 , \ov{e}_2) : \ov{M}_{0;2,0}( \X ,1) \to C \times C \cong \PP^1 \times \PP^1$$ is an isomorphism (coarsely).  A point $(\tilde{f}: \tilde{\Sigma} \to F, \Sigma, R_1, R_2) \in \ov{M}_{0;2,0}( \X ,1)$ parameterizes a curve of minimal degree, so at most one component $\Sigma_1$ of the base curve $\Sigma$ is not collapsed.  Since there are only two marked points, either $\Sigma=\Sigma_1$ or $\Sigma=\Sigma_1 \amalg_N \Sigma_2$ with $R_1,R_2 \in \Sigma_2 \setminus \{ N \}$ (stability!).  

In the first case, $R_1,R_2 \in \Sigma \cong \PP^1$ so $\tilde{\Sigma} \cong \PP^1$ and $\tilde{f} : \tilde{\Sigma} \to F$ is equivariant of bidegree $(1,1)$.  In particular, $\tilde{f}$ is an embedding so $$\ov{e}=((A \subset A^\perp),(B \subset B^\perp)) \in C \times C$$ with $A \neq B$.  Now, notice that any equivariant map $g:\PP^1 \to F$ is determined by $\pi_1 g : \PP^1 \to \PP^2$, so since $\pi_1 \tilde{f}$ above is of degree $1$, it is determined by $A,B \in \PP^2$ hence there is at most one map in $\ov{M}_{0;2,0}( \X ,1)$ with no collapsed components with $\ov{e}$ as above.  To see that there is exactly one, let $V := {\rm Span}(A,B)$ so $\CC^3=V \oplus V^\perp$ and consider the map $\PP V \to F$ given by $$C \mapsto (C \subset C \oplus V^\perp) \in F.$$  This is an embedding with $\ZZ_2$-invariant image because $$((C \oplus V^\perp)^\perp \subset C^\perp) = (D \subset D \oplus V^\perp)$$ for some $D \subset V$.  Notice that this curve is of bidegree $(1,1)$ and meets $C$ exactly at the two coordinates of $\ov{e}$.  

When there is a collapsed component $\Sigma_2 \cong \PP^1$, it contains the two marked points so $$\ov{e} = ((A \subset A^\perp),(A \subset A^\perp)) \in C \times C$$ and the cover $\tilde{\Sigma}$ is unramified over the uncollapsed component $\Sigma_1$, so over $\Sigma_1$, we have $\tilde{f}|_{\tilde{\Sigma}_1}: \PP^1 \amalg \PP^1  \to F$ with degrees $(0,1)$, $(1,0)$.  Since $\tilde{f}$ is equivariant, it is uniquely determined by, say, the degree $(0,1)$ map, which must be an isomorphism onto $$\{ (A \subset B) \} \subset F$$ ($A$ fixed, $B$ varying).  The degree $(1,0)$ map must be an isomorphism onto $$\{ (C \subset A^\perp) \} \subset F,$$ so the two preimages $N_1, N_2 \in \tilde{\Sigma}$ of the node $N \in \Sigma$ map to $(A \subset A^\perp)$, and $\tilde{\Sigma_2}$ is a $\PP^1$ glued to $\tilde{\Sigma_1}$ at $N_1, N_2$ and collapsed by $\tilde{f}$ to $(A \subset A^\perp)$.  This proves the claim, so by the Divisor Axiom (see Page 193 in \cite{CK} for the axioms of Gromov-Witten theory used in this paper) we compute $$\langle S_1,S_4,S_4 \rangle^1 = (1/2) \langle p_1+p_2,W_1+W_2 \rangle \int_{C \times C} \pi_1^*x \cdot \pi_2^*x  = 1.$$

For dimension reasons, only nine 3-point Gromov-Witten numbers could possibly be non-zero.  They turn out to be: $$\begin{array}{cc} {\rm Nontwisted \,} (r=0) & {\rm Twisted \,} (r=2) \\ \langle S_1,S_3,S_3 \rangle^1 = 9 & \langle S_1,S_4,S_4 \rangle^1 = 1 \\ \langle S_1,S_1,S_5 \rangle^1 = 6 & \langle S_2,S_3,S_4 \rangle^1 = 3 \\ \langle S_3,S_3,S_5 \rangle^2 = 54 & \langle S_2,S_2,S_5 \rangle^1 = 6 \\ \langle S_1,S_5,S_5 \rangle^2 = 36 & \langle S_4,S_4,S_5 \rangle^2 = 3 \\ \langle S_5,S_5,S_5 \rangle^3 = 0 \end{array}$$ We computed 3 of these above, and the rest can be derived from associativity of the small quantum orbifold cohomology ring, by considering the following associativity checks in order: $$\begin{array}{cccccrcl} (S_1 \star S_1) \star S_2 & = & S_1 \star (S_1 \star S_2) & \Longrightarrow & \langle S_2,S_3,S_4 \rangle^1 & = & 3 \\ (S_2 \star S_3) \star S_1 & = & S_2 \star (S_3 \star S_1) &\Longrightarrow & \langle S_2,S_2,S_5 \rangle^1 & = & 6 \\ (S_3 \star S_2) \star S_2 & = & S_3 \star (S_2 \star S_2) & \Longrightarrow & \langle S_3,S_3,S_5 \rangle^2 & = & 54 \\ (S_1 \star S_3) \star S_5 & = & S_1 \star (S_3 \star S_5) & \Longrightarrow & \langle S_1,S_5,S_5 \rangle^2 & = & 36 \\ & & & \Longrightarrow & \langle S_5,S_5,S_5 \rangle^3 & = & 0 \\ (S_4 \star S_1) \star S_3 & = & S_4 \star (S_1 \star S_3) & \Longrightarrow & \langle S_4,S_4,S_5 \rangle^2 & = & 3 \end{array}$$  The same result can be obtained from the orbifold WDVV equation\footnote{It is also interesting to compute the degree $1$ twisted invariants by localization.}.  These 3-point invariants form, in the usual way, a commutative (associative!) graded ring structure on $H^*_{orb}( \X ) \otimes_{\QQ} \QQ[q]$ whose multiplication table is given below.  $$\begin{array}{lll} S_1 \star S_1 = S_3+2q & \quad & S_2 \star S_2 = (2/3)S_3+2q \\ S_1 \star S_2 = 4S_4 & \quad & S_2 \star S_3 = 6qS_2 \\ S_1 \star S_3 = S_5+3qS_1 & \quad & S_2 \star S_4 = (1/6)S_5+qS_1 \\ S_1 \star S_4 = 2qS_2 & \quad & S_2 \star S_5 = 12qS_4 \\ S_1 \star S_5 = 2qS_3+12q^2 & & \\ S_3 \star S_3 = 3qS_3+18q^2 & \quad & S_4 \star S_4 = (1/3)qS_3 + q^2 \\ S_3 \star S_4 = 6qS_4 & \quad & S_4 \star S_5 = 6q^2S_2 \\ S_3 \star S_5 = 18q^2S_1 & \quad & S_5 \star S_5 = 12q^2S_3\end{array}$$   Deforming the relations in $H^*_{orb}( \X )$ we can give a presentation: $$QH^*_{orb}( \X ) = \QQ[S_1,S_2,q]/ \langle S_2^3-6qS_2,3S_2^2-2S_1^2-2q \rangle$$

\end{section}

\begin{section}{The Crepant Resolution $Y$} \label{section:y}

After fixing a non-degenerate symmetric bilinear form $\langle,\rangle:\CC^3 \otimes \CC^3 \to \CC$ we may regard the smooth variety $F$ of complete flags in $\CC^3$ as the subspace of $\PP^2 \times \PP^2$ consisting of pairs $(A,B)$ with $\langle A,B \rangle=0$.  The $\ZZ_2$-action interchanging the factors of $\PP^2 \times \PP^2$ restricts to the $\ZZ_2$-action on $F$ taking a flag to its orthogonal complement.  The rational map $\PP^2 \times \PP^2 \dashrightarrow \PP^{2*}$ taking two (distinct) 1-dimensional subspaces of $\CC^3$ to the 2-dimensional subspace they span is undefined on the diagonal, but can be resolved by blowing up the diagonal (the $\ZZ_2$-action lifts to an action on the blowup).  The fiber of the resulting map $p: {\rm Bl}_{\Delta} (\PP^2 \times \PP^2) \to \PP^{2*}$ over $W \in \PP^{2*}$ is canonically $\PP W \times \PP W$, while the fiber $p^{-1}(W) \cap \tilde{F}$ of the restriction of $p$ to the proper transform $\tilde{F}$ of $F$ is $$\{ (A,B) \in \PP W \times \PP W : \langle A,B \rangle = 0 \}.$$  The topology of this fiber depends on whether $W^\perp \subset W$.  If so, then the fiber is just $$\{ (W^\perp, A) : A \subset W \} \cup_{ (W^\perp,W^\perp) } \{ (A,W^\perp) : A \subset W \} \cong \PP^1 \lor \PP^1,$$ while if $W^\perp \cap W = (0)$, then the fiber is the graph of the idempotent automorphism of $\PP W$ taking $A$ to $A^\perp \cap W$.  The map $p$ is $\ZZ_2$-equivariant, so it descends to a map on the $\ZZ_2$-quotients.  The fiber of $\pi$ over $W$ is of course $\Sym^2 \PP W \cong \PP \Sym^2 W$.  T.~Graber~\cite{G} used this to note that ${\rm Hilb}^2 \PP^2 = {\rm Bl}_{\Delta} (\PP^2 \times \PP^2) / \ZZ_2$ is the projectivization of $\Sym^2 W$, where $W$ now denotes the tautological rank 2 bundle on $\PP^{2*}$.  Since the first Chern class of $\OO^3/W$ is the positive generator $T_1$ of $H^2(\PP^{2*},\ZZ) \cong \ZZ$, the SES $$0 \to W \to \OO^3 \to \OO^3/W \to 0$$ shows $c(W)=1-T_1+T_1^2$ from which one easily computes $c(\Sym^2 W)=1-3T_1+6T_1^2$.  We will argue in a moment that $Y := \tilde{F} / \ZZ_2$ is a crepant resolution of singularities of $F / \ZZ_2$.  For now, notice that the fiber of $\pi|_Y$ over $W$ is canonically the projectivization of the 2-dimensional space $$V := \{  w_1 \lor w_2 \in \Sym^2 W :  \langle w_1,w_2 \rangle = 0 \} \subset \Sym^2 W$$ so that we may describe $Y$ as the projectivization of the corresponding rank 2 vector bundle $V \subset \Sym^2 W$ on $\PP^{2*}$.  The quotient bundle $\Sym^2 W / V$ is isomorphic to the trivial bundle by the map $[w_1 \lor w_2] \mapsto \langle w_1,w_2 \rangle$ so we have $c(V)=c(\Sym^2 W)=1-3T_1+6T_1^2$.

The proper transform $\tilde{F}$ of $F \subset \PP^2 \times \PP^2$ in ${\rm Bl}_\Delta (\PP^2 \times \PP^2)$ is obtained by blowing up $F$ along $K=(F \cap \Delta) \cong \PP^1$, which is the locus of flags of the form $(W^\perp \subset W)$ (i.e. the fixed locus of the $\ZZ_2$-action on $F$).  In local analytic coordinates, this blowup is just $({\rm Bl}_{(0,0)} \AA^2) \times \AA^1$.  Then we take the $\ZZ_2$-quotient to obtain $Y$.  One easily checks locally that blowing up the $\ZZ_2$-quotient of $F$ along $K$ yields a crepant resolution (locally $F / \ZZ_2$ looks like an $A_1$ singularity times $\AA^1$).  However, the order in which we do the blowup and take the $\ZZ_2$ quotient doesn't matter.  Locally this corresponds to the (easily checked) commutativity of the natural diagram $$\begin{CD} {\rm Tot \;} \OO(-1) = {\rm Bl}_{(0,0)} \AA^2 @>>> \AA^2 \\ @VVV @VVV \\ {\rm Tot \;} \OO(-2) = {\rm Bl}_{(0,0)} \AA^2 / \ZZ_2 @>>> \AA^2 / \ZZ_2 \end{CD}$$ where the vertical maps are $\ZZ_2$-quotients and the horizontal maps are blow-ups.

Applying the Leray-Hirsch theorem to the projective bundle description of $Y$ yields a presentation of its cohomology ring $$H^*(Y,\ZZ) = \ZZ[T_1,T_2]/\langle T_1^3,T_2^2-3T_1T_2+6T_1^2 \rangle$$ where $T_1$ is pulled back from $\PP^{2*}$ and $T_2 = c_1( \OO_{\PP V}(1) )$.  We can also compute the Chern classes $c(TY) = 1+2T_2-6T_1^2+6T_1T_2+6T_1^2T_2$.

\bigskip

\centerline {\bf Curve Classes in $Y$}

\bigskip

Restricting $\PP V$ to a line in $\PP^{2*}$ yields a Hirzebruch surface whose algebraic type turns out to depend on the type of line as follows.  A \emph{generic} line is of the form  $$\ell_A := \{ B \in \PP^{2*} : A \subset B \}$$ for some fixed 1-dimensional $A \subset \CC^3$ where $A \cap A^\perp = (0)$.  Consider the line bundle $L \subset W|_{\ell_A}$ over $\ell_A$ whose fiber over $B \in \ell_A$ is $B \cap A^\perp$.  Choose some nonzero $a \in A$ and notice that $1 \mapsto a$ gives an injection $\OO \hookrightarrow W|_{\ell_A}$, yielding a splitting $W|_{\ell_A} \cong \OO \oplus L$, so $c_1(L)=-1$.  There is also an injective map of vector bundles on $\ell_A$ $$\begin{array}{c} L \to V|_{\ell_A} \\ b \mapsto a \lor b\end{array}$$ and $V|_{\ell_A}$ has first Chern class $-3$ so the quotient line bundle has first Chern class $-2$, hence the corresponding SES is split for cohomological reasons and we get $$V|_{\ell_A} \cong \OO(-1) \oplus \OO(-2).$$  We conclude that the projectivization of $V|_{\ell_A}$ is the Hirzebruch surface $\FF_1$.  

Now consider some $A \in \PP^2$ with $A \subset A^\perp$ (equivalently $\langle A,A \rangle = 0$).  We call the corresponding line $\ell_A \subset \PP^{2*}$ a \emph{jump} line.  Here we may simply take a nonzero $a \in A$ and get an injective vector bundle morphism $\OO \to V|_{\ell_A}$ by taking $1 \to a \lor a$.  The quotient has first Chern class $-3$ so the SES is again split for cohomological reasons and we have $$V|_{\ell_A} \cong \OO \oplus \OO(-3)$$ so $\PP V|_{\ell_A}$ is the Hirzebruch surface $\FF_3$.  On jump lines, the dimensions of $H^0(\ell_A,V|_{\ell_A})$ and $H^1(\ell_A, V|_{\ell_A})$ jump up from $0,1$ (for a generic line) to $1,2$ respectively.

We now study rational curves in $Y$.  Identify a curve class $[C]$ with the pair $(\langle T_1,[C] \rangle, \langle T_2,[C] \rangle) \in \ZZ^2$.  Using WDVV (see \cite{CK} or \cite{AGV2} for orbifolds) we will see that the 3-point Gromov-Witten numbers of $Y$ can be computed by studying only the moduli spaces $\ov{M}_{0;2}(Y,(0,1))$, $\ov{M}_{0;3}(Y,(0,2))$, and $\ov{ M}_{0;0}(Y,(n,0))$.  Curves corresponding to pairs $(0,a)$ are collapsed by $\pi$ and are thus (branched) covers of some fiber $F$ of $\pi$.  Since $\langle T_2,[F] \rangle = 1$ we always have $a \geq 0$ and the moduli space $$\ov{ M}_{0;n}(Y,(0,a))$$ is a fiber bundle over $\PP^{2*}$ whose fiber over $W$ is $\ov{M}_{0;n}(\PP W, a)$.  This is smooth (as a stack) of the expected dimension $$\dim Y - 3 + n + \langle c_1(TY) , (0,a) \rangle = n - 2a.$$  In particular, $$\ov{M}_{0;2}(\PP W, 1)=\PP W \times \PP W$$ so we have a simple description of $\ov{ M}_{0;2}(Y,(0,1))$ and its evaluation maps, from which we can easily evaluate the 2-point Gromov-Witten invariants for the homology class $(0,1)$.  The only such invariants which are non-zero are $$\begin{array}{ccc} \langle T_2,T_1^2T_2 \rangle^{0,1} = 1 & \quad {\rm and} \quad \, &  \langle T_1T_2, T_1T_2 \rangle^{0,1} = 1. \end{array}$$  Using the Divisor Axiom, we can also evaluate all 3-point invariants for this homology class.  

Since $\ov{ M}_{0;3}(Y,(0,2))$ is a fiber bundle over $\PP^{2*}$ we know that the 3-point invariants of the form $\langle T_1^iT_2^l, T_1^jT_2^m, T_1^kT_2^n \rangle^{0,2}$ will vanish if $i+j+k > 2$.  Combining this with the Dimension Axiom shows that all such 3-point invariants are zero.

Now we turn our attention to curves contained in one of the Hirzebruch surfaces mentioned above.  The effective cone of a Hirzebruch surface is generated by the fiber class and the class of the \emph{rigid section}: any Hirzebruch surface can be written as $\PP (\OO \oplus \OO(n))$ where $n \leq 0$, and the rigid section $s$ is obtained by taking the subspace spanned by the trivial factor.  As long as $n \neq 0$ this is the unique section in its homology class.  The rigid section has normal bundle $\OO(n)$.  One can easily check that the rational curve corresponding to the rigid section $s$ of the Hirzebruch surface $\FF_1$ over a generic line $\ell_A$ corresponds to the homology class $(1,1)$, while the rational curve corresponding to the rigid section $t$ of $\FF_3$ over a jump line $\ell_A$ is in the homology class $(1,0)$.  

Now consider the moduli space $\ov{ M}_{0;0}(Y,(1,0))$.  The rigid section $t$ of $\FF_3$ is unique, so this moduli space is the space of jump lines, which is a $\PP^1$ given by the conic $S_2 \subset \PP^2$ consisting of those $A \in \PP^2$ with $\langle A,A \rangle = 0$.  We may explicitly identify $M_1 := \ov{ M}_{0;1}(Y,(1,0))$ as well.  The map $\pi_1$ forgetting the marked point makes $M_1$ a $\PP^1$-bundle over $\ov{ M}_{0;0}(Y,(1,0)) \cong S_2 \cong \PP^1$.  I claim this is a trivial Hirzebruch surface.  First of all, for $A \in S_2$, giving a point of the section $t_A$ is the same as giving a point of the line $\ell_A$, which is the same as giving a 1-dimensional subspace of $\CC^3 / A$.  However, it is a simple matter to check that there is a canonical isomorphism $\PP( \CC^3/A ) = \PP \Hom(A, \CC^3 / A)$, so we need only show that the restriction of the tangent bundle of $\PP^2$ to $S_2$ has balanced splitting type.  Indeed, this is a special case of the well-known fact that the restriction of $T\PP^n$ to any rational normal curve (image of a degree $n$ embedding $\PP^1 \to \PP^n$) has balanced splitting type.

The expected dimension of the moduli space $\ov{ M}_{0;0}(Y,(1,0))$ is $0$ (this is a general phenomenon for crepant resolutions---see below), so we must identify the virtual fundamental class.  This is given by the (Poincar\'e dual of the) first Chern class of the vector bundle over $S_2$ whose fiber over $A$ is $H^1( t_A , N_{ t_A / Y })$.  We will soon check that the rank of this vector bundle is 1 (the excess dimension).  The moduli space $M_1$ is contained in $Y$ as the union of all rigid sections over jump lines (these are disjoint).  Since $M_1$ is a trivial Hirzebruch surface, the normal bundle of a rigid section $t_A$ in $M_1$ is trivial.  The SES  $$0 \to Tt_A \to TY|_{t_A} \to N_{t_A / Y} \to 0$$ implies that $c_1(N_{t_A / Y})=-2$.  Thus the SES $$0 \to N_{t_A / M_1} \cong \OO_{t_A} \to N_{t_A / Y} \to N_{M_1 / Y}|_{t_A} \to 0$$ must be split so that $N_{t_A / Y} \cong \OO_{t_A} \oplus \OO_{t_A}(-2)$ and we have a canonical isomorphism $$H^1(t_A,N_{t_A / Y}) = H^1(t_A, N_{M_1 / Y}|_{t_A}).$$

We introduce an algebraic $\CC^*$-action on $Y$ by taking a maximal torus in the Lie group ${\rm SO}_3(V,\langle , \rangle)$ of matrices preserving $\langle,\rangle$, which acts naturally on $Y$.  Explicitly, we may take the $\CC^*$-action $$\lambda \cdot (z_0,z_1,z_2) := (\lambda z_0, z_1, \lambda^{-1} z_2)$$ on $\CC^3$, which preserves the bilinear form $$\langle (z_0,z_1,z_2),(w_0,w_1,w_2) \rangle := z_0w_2+z_1w_1+z_2w_0$$ mentioned in the introduction.  The jump lines fixed by this action correspond to the points $[1:0:0] \in S_2$ and $[0:0:1] \in S_2$.  This $\CC^*$-action on $Y$ has $6$ fixed points.  The cohomology classes $T_1,T_2,$ and $c_1(N_{M_1 / Y})$ naturally lift to equivariant cohomology classes ($M_1$ is invariant under the action so the last class lies in $H^*_{\CC^*}(M_1)$) whose weights at the fixed points are listed in the table below.

$$\begin{array}{ccccc} P \in Y^{\CC^*} & T_P Y & T_1|_P & T_2|_P & N_{M_1 / Y}|_P \\
P_1 := \langle e_0 \lor e_0 \rangle \subset W|_{E_{01}} \subset \Sym^2 E_{01} & -1,-2,-1 & -1 & -2 & -2 \\
P_2 := \langle e_0 \lor e_1 \rangle \subset W|_{E_{01}} \subset \Sym^2 E_{01} & -1,-2,1 & -1 & -1 &  \\
P_3 := \langle e_2 \lor e_2 \rangle \subset W|_{E_{02}} \subset \Sym^2 E_{02} & 1,-1,4 & 0 & 2 & 4 \\
P_4 := \langle e_0 \lor e_0 \rangle \subset W|_{E_{02}} \subset \Sym^2 E_{02} & 1,-1,-4 & 0 & -2 & -4 \\
P_5 := \langle e_1 \lor e_2 \rangle \subset W|_{E_{12}} \subset \Sym^2 E_{12} & 1,2,-3 & 1 & 1 &  \\
P_6 := \langle e_2 \lor e_2 \rangle \subset W|_{E_{12}} \subset \Sym^2 E_{12} & 1,2,1 & 1 & 2 & 2 \end{array}$$

We will use this to show that $N_{M_1 \subset Y} \cong \OO_{\PP^1 \times \PP^1}(6,-2)$.  The map $$[s:t] \mapsto \langle e_0 \lor e_0 \rangle \subset W|_{ {\rm Span} (e_0,(0,s,t)) } \in M_1$$ is the inclusion of the fiber of $\pi_1$ over $E_0 \in S_2 \subset \PP^2$.  This fiber contains two fixed points: $P_1$ and $P_4$.  From the localization chart we can see that $$\langle c_1(N_{M_1 / Y}) , [ \pi_1^{-1}(E_0) ] \rangle = -2$$ so $N_{M_1 / Y} = \OO_{\PP^1 \times \PP^1 }(d,-2)$ for some integer $d$.  The map $$[s:t] \mapsto \langle (s^2 , \sqrt{-2}st , t^2) \lor (s^2 , \sqrt{-2}st , t^2) \rangle \subset W|_{ \Ker (z \mapsto \langle z, (s^2 , \sqrt{-2}st , t^2) \rangle ) } \in M_1 $$ is a $\CC^*$-invariant section of $\pi_1$ containing the fixed points $P_1$ and $P_3$, so it is either a $\CC^*$-invariant curve in $M_1$ of degree $(1,1)$ or it is a fiber of $\pi_2$.  The integral of $c_1(N_{M_1 / Y})$ over this curve is $6$, so in the first case we would have $N_{M_1 \subset Y} \cong \OO_{\PP^1 \times \PP^1}(8,-2)$ and in the second case we would have $N_{M_1 \subset Y} \cong \OO_{\PP^1 \times \PP^1}(6,-2)$.  However, the first case is impossible because then the fiber of $\pi_2$ containing $P_1$ would also contain $P_6$ and thus the integral over this fiber would be $4$ (not $8$), whereas, in the second case, the curve containing $P_1$ and $P_6$ is of degree $(1,1)$ and the integral of $c_1(N_{M_1 / Y})$ is correctly given by $4$.

\newpage

\centerline {\bf The Resolution Map $r:Y \to F/\ZZ_2$}

\bigskip

We may explicitly describe the crepant resolution map $r:Y \to F/\ZZ_2$ as follows.  A point of $Y$ is specified by some $W \in \PP^{2*}$ together with a subspace $A=\langle a \lor b \rangle \subset V|_W$.  The image of this point under $r$ will be $$[\langle a \rangle \subset \langle b \rangle^\perp ] = [ \langle b \rangle \subset \langle a \rangle^\perp ] \in F / \ZZ_2.$$  It is easy to check that this is well-defined.  To see that this is an isomorphism on the locus of $[A_1=\langle a \rangle \subset A_2] \in F / \ZZ_2$ where $A_2 \neq A_1^\perp$, just notice that $A_1^\perp \cap A_2$ is a 1-dimensional subspace of $\CC^3$, spanned, say, by $b$ so that $$r^{-1}([A_1=\langle a \rangle \subset A_2]) = \langle a \lor b \rangle \in V|_{\langle a,b \rangle}.$$

Recall that we constructed the rigid section over a jump line $\ell_A$ by always taking the 1-dimensional subspace $[a \lor a]$ in the fiber and letting only the 2-dimensional subspace $W$ vary (over all 2-dimensional subspaces containing $A$).  Thus the resolution map collapses the rigid section over a jump line to a point so that the corresponding homology class $(1,0)$ is on the boundary of the effective cone of $Y$.  Since the projection to $\PP^{2*}$ collapses the fiber class, it also lies on the boundary of the effective cone, thus $H_2^{\rm eff}(Y,\ZZ) = \{ (a,b) : a,b \geq 0 \}.$

Since $r$ is a crepant resolution (i.e. $\omega_Y = r^* \omega_{F / \ZZ_2}$) we expect to have a degree $0$ quantum parameter because the canonical bundle of $Y$ will evaluate $0$ on the curve class collapsed by $r$ (here: the class $(1,0)$).  The same phenomenon occurs for the crepant resolution ${\rm Hilb}^2 \PP^2 \to \Sym^2 \PP^2$ where the class of $0$-dimensional subschemes supported at a fixed point is collapsed.

\bigskip

\centerline {\bf Classical Geometry of the Quotient}

\bigskip

Here we make contact with  some classical geometry by giving another description of the quotient $F/ \ZZ_2$ and the resolution $r: Y \to F / \ZZ_2$; none of this is strictly necessary in what follows.  Recall that $F$ is $\ZZ_2$-equivariantly embedded in $\PP V \times \PP V$ (with $\ZZ_2$ exchanging the factors) as the set $\{ (A,B) : \langle A,B \rangle = 0 \}$.  Then $$(\PP V \times \PP V) / \ZZ_2 \hookrightarrow \PP \Sym^2 V$$ via the Segre embedding.  The image of $F / \ZZ_2$ under this embedding inside the image of $(\PP V \times \PP V) / \ZZ_2$ is given by a single linear equation.  For example, if $V = \CC^3$ with the ``back-to-front" inner product from the introduction, then the Segre embedding is $$[a:b:c]+[d:e:f] \mapsto [ad:ae+bd:be:af+cd:cf:bf+ce]$$ and the image of $F / \ZZ_2$ is the intersection of the image of $( \PP^2 \times \PP^2 ) / \ZZ_2$ with the hyperplane $H=V(X_2 + X_3) \cong \PP^4$.  The fixed locus $C \subset F$ is a curve in $\PP^2 \times \PP^2$ of bidegree $(2,2)$ so its image (or rather, the image of $C / \ZZ_2 = C \subset F / \ZZ_2$) under the Segre embedding is a rational normal curve in $H$.  In the above coordinates, $$C = \{ [2s^2:2ist:t^2]+[2s^2:2ist:t^2] \}$$ and the rational normal curve is $$\{ [4s^4 : 8is^3t : -4s^2t^2 : 4s^2t^2 : t^4 : 4ist^3 ] \}.$$ We will show that the image of $F / \ZZ_2$ is the secant (or \emph{chordal}) variety of this rational normal curve in $H$, which is a singular degree $3$ hypersurface (see page 120 in \cite{JH}).  Indeed, taking the $SO_3$-action into account, it is enough to show that the image of $F / \ZZ_2$ contains the line between, say, the images of $[1:0:0]+[1:0:0]$ and $[0:0:1]+[0:0:1]$, as well as the tangent line to the rational normal curve at, say, the image of $[1:0:0]+[1:0:0]$.  This first line is $\{ [s:0:0:0:t:0] \}$, which is the image of $$\{ [ \sqrt{s} : 0 : \sqrt{-t}] + [ \sqrt{s} : 0 : - \sqrt {-t} ] \} \subset F / \ZZ_2.$$  The second line is $\{ [s:t:0:0:0:0 ] \}$, which is the image of $\{ [1:0:0]+[s:t:0] \}.$

Every line in the chordal variety $F / \ZZ_2$ is either (1) a tangent line to the rational normal curve or (2) a line connecting two distinct points on the rational normal curve.  Taking the $SO_3$-action into account and condsidering the explicit computations above we see that the preimage of a line of type (1) in $F$ consists of two smooth rational curves of degree $(1,0)$ and $(0,1)$ exchanged by the $\ZZ_2$-action and meeting at a point of $C$; the preimage of a line of type (2) in $F$ is a $\ZZ_2$-invariant smooth rational curve of degree $(1,1)$ meeting $C$ at two distinct points.

\end{section}

\begin{section}{Gromov-Witten Theory of $Y$} \label{section:gwy}

In this section, we compute the 3-point Gromov-Witten invariants of $Y$ via WDVV.  Using the basis $\{ T_i \}$ and its Poincar\' e dual basis $\{ T^i \}$ below $$\begin{array}{c} T_0=1 \\ T_1 \quad T_2  \\ T_3 = T_1^2 \quad T_4 = T_1T_2 \\ T_5 = T_1^2T_2 \end{array} \quad \begin{array}{c} T^0=T_5 \\ T^1=T_4-3T_3 \quad T^2 = T_3 \\ T^3 = T_2-3T_1 \quad T^4 = T_1 \\ T^5 = T_0 \end{array}$$ for the cohomology of $Y$, the Poincar\'e duality metric $G=(G_{ij})$ and its inverse $G^{-1} = (G^{ij})$ are given as below.  $$\begin{array}{lcr} G = \begin{pmatrix} 0&0&0&0&0&1 \\ 0&0&0&0&1&0 \\ 0&0&0&1&3&0 \\ 0&0&1&0&0&0 \\ 0&1&3&0&0&0 \\ 1&0&0&0&0&0 \end{pmatrix} & \quad & G^{-1} = \begin{pmatrix} 0&0&0&0&0&1 \\ 0&0&0&-3&1&0 \\ 0&0&0&1&0&0 \\ 0&-3&1&0&0&0 \\ 0&1&0&0&0&0 \\ 1&0&0&0&0&0 \end{pmatrix} \end{array}$$ WDVV says that for any $\phi_1,\phi_2,\phi_3,\phi_4 \in H^*(Y)$, any $n \geq 0$, any $\gamma_1,\dots,\gamma_n \in H^*(Y)$, and any $2$-dimensional homology class $\beta$ we have \begin{eqnarray} \nonumber & & \sum_{ \substack{ a,b \\ \beta_1 + \beta_2 = \beta \\ V \cup W = [n] }} \langle \phi_1,\phi_2,T_a,\gamma_{v_1},\dots,\gamma_{v_k} \rangle^{\beta_1}G^{ab} \langle \phi_3,\phi_4,T_b,\gamma_{w_1},\dots,\gamma_{w_{n-k}} \rangle^{\beta_2}  \\  \label{wdvv} & = & \sum_{ \substack{ a,b \\ \beta_1 + \beta_2 = \beta \\ V \cup W = [n] }} \langle \phi_1,\phi_4,T_a,\gamma_{v_1},\dots,\gamma_{v_k} \rangle^{\beta_1}G^{ab} \langle \phi_2,\phi_3,T_b,\gamma_{w_1},\dots,\gamma_{w_{n-k}} \rangle^{\beta_2}.\end{eqnarray}  The terms on the LHS where either $\beta_1$ or $\beta_2$ is zero sum to give \begin{eqnarray} \label{lhs} \LHS^\beta & := & \langle \phi_1,\phi_2,\phi_3 \phi_4,\gamma_1,\dots,\gamma_n \rangle^\beta + \langle \phi_3, \phi_4, \phi_1 \phi_2, \gamma_1,\dots,\gamma_n \rangle^\beta \end{eqnarray} and similarly the terms on the RHS with $\beta_1$ or $\beta_2$ equal to zero give \begin{eqnarray} \label{rhs} \RHS^\beta := \langle \phi_1,\phi_4,\phi_2 \phi_3,\gamma_1,\dots,\gamma_n \rangle^\beta + \langle \phi_2, \phi_3, \phi_1 \phi_4 ,\gamma_1,\dots,\gamma_n \rangle^\beta .\end{eqnarray}  In the remainder of this section we will always use this equation with $n=0$ and with the insertions $\phi_1,\phi_2,\phi_3,\phi_4$ equal to some $T_i,T_j,T_k,T_l$, so we will specify a WDVV equation by indicating the choice of $\beta$ and $(i,j,k,l)$.

We begin by computing all Gromov-Witten invariants for homology classes of the form $(n,0)$.  By the Divisor and Dimension Axioms we need only compute the $0$-point invariant $\langle \rangle^{n,0}$.  Since $(1,0)$ is collapsed by the resolution map $r$, the curves of class $(n,0)$ are just $n$-fold branched covers of rigid sections over jump lines.  Such a map is specified by a point in the first factor of $$M_1 = \ov{M}_{0;0}(Y,(1,0)) \cong \PP^1 \times \PP^1$$ together with an element of $\ov{M}_{0;0}(\PP^1,n)$ so that the moduli space $\ov{M}_{0;0}(Y,(n,0))$ (which has expected dimension $0$) is a product $\PP^1 \times \ov{M}_{0;0}(\PP^1,n)$.  Furthermore, we can identify the obstruction class as the Euler class of the vector bundle $$\pi_1^* \OO_{\PP^1}(6) \otimes \pi_2^* H^1(C,f^* \OO_{\PP^1}(-2))$$ where, by abuse of notation, the second factor denotes the vector bundle on $\ov{M}_{0;0}(\PP^1,n)$ whose fiber over $(f,C)$ is $H^1(C,f^* \OO_{\PP^1}(-2))$.  Recall that for a $2$-dimensional vector space $V$ we have a natural SES $$0 \to \OO_{\PP V}(-1) \to V \to V / \OO_{\PP V}(-1) \to 0$$ which we can twist by $\OO_{\PP V} (-1)$ and pull-back by any map $f$ to get a SES $$0 \to f^* \OO_{\PP V}(-2) \to  V \otimes f^* \OO_{\PP V} (-1) \to f^* [ (V / \OO_{\PP V} (-1)) \otimes \OO_{\PP V} (-1)] \to 0$$ where the vector bundle on the right is (non-canonically) trivial.  The associated LES in cohomology gives $$0 \to H^0(C,\OO_C) \to H^1(C,f^* \OO_{\PP V} (-2)) \to H^1(C,V \otimes f^* \OO_{\PP V}(-1)) \to 0.$$  Thus we compute $$\begin{array}{rcl} \langle \rangle^{n,0} &=& \int_{ \ov{M}_{0;0}(Y,(n,0))} e( {\rm Obstruction \; Bundle} ) \\ 
&=& \int_{ \PP^1 \times \ov{M}_{0;0}(\PP^1,n) } c_{2n-1}( \pi_1^* \OO_{\PP^1}(6) \otimes \pi_2^* H^1(C,f^* \OO_{\PP^1} (-2) ))\\ 
&=& 6 \int_{ \ov{M}_{0;0}(\PP^1,n)} c_{2n-2} ( H^1(C,f^* \OO(-1) \oplus \OO(-1)) ) \\ &=& 6/n^3 \end{array}$$ where the last equality is the Aspinwall-Morrison formula (Theorem~9.2.3 in \cite{CK}).  Thus the only non-zero $3$-point invariant for the homology class $(n,0)$ is $$\langle T_1, T_1, T_1 \rangle^{n,0} = 6.$$

Next we compute the $3$-point invariants when $\beta = (n,1)$.  Using the Dimension and Divisor Axioms we see that it suffices to compute the $4$ numbers $$\begin{array}{cc} a_{225}^n := \langle T_2,T_2,T_5 \rangle^{n,1} &  a_{233}^n := \langle T_2,T_3,T_3 \rangle^{n,1} \\
 a_{234}^n := \langle T_2,T_3,T_4 \rangle^{n,1} & a_{244}^n := \langle T_2,T_4,T_4 \rangle^{n,1} \end{array}$$ because the only other such $3$-point invariants that might be non-zero are then determined by the Divisor Axiom.  For example $$\langle T_1,T_1,T_5 \rangle^{n,1} = n^2 \langle T_5 \rangle^{n,1} = n^2 \langle T_2,T_2,T_5 \rangle^{n,1} = n^2 a_{225}^n $$ and similarly: $$\begin{array}{cc}  \langle T_1,T_2,T_5 \rangle^{n,1} = na_{225}^n &  \langle T_1,T_3,T_3 \rangle^{n,1} = na_{233}^n \\  \langle T_1,T_3,T_4 \rangle^{n,1} = na_{234}^n &  \langle T_1,T_4,T_4 \rangle^{n,1} = na_{244}^n \end{array}$$   We computed these invariants for $n=0$ in Section~\ref{section:y}.  We can get a system of four equations yielding a recursive formula (in $n$) for these $4$ invariants by applying WDVV with $\beta=(n,1)$ and $ (i,j,k,l) =(1,1,2,3), \, (1,1,2,4), \, (1,2,2,3), \, (1,2,2,4).$  Writing each of these equations in the form $$\LHS^\beta - \RHS^\beta = C$$ where $C$ is determined by $3$-point invariants for homology classes $\alpha < (n,1)$ (in at least one coordinate) we get a system of $4$ equations $$\begin{array}{cccccc} n^2a^n_{225} & +a^n_{233} & -na^n_{234} & & = & C^n_1 \\ (3n^2-n)a^n_{225} & & +a^n_{234} & -na^n_{244} & = & C_2^n \\ na^n_{225} & +6na^n_{233} & +(1-3n)a^n_{234} & & = & 0 \\ (3n-1)a^n_{225} & & +6na^n_{234} & + (1-3n)a^n_{244} & = & 0 \end{array}$$ which uniquely determines the 4 unknowns in terms of $C^n_1$ and $C^n_2$ when $n > 0$ because the determinant of the coefficient matrix is $n^2(6n-1)(3n^2-6n+1)$.  The coefficients $C_1^n$ and $C_2^n$ are easily worked out because if $\beta_1 + \beta_2 = (n,1)$, then one of $\beta_1$ or $\beta_2$ is of the form $(d,0)$ and the corresponding three point invariant is almost always $0$.  We get $$\begin{array}{ccc} C^n_1  =  6 \sum_{d=0}^{n-1} (3a^d_{233} - a^d_{234}) & \quad &
C^n_2  =  6 \sum_{d=0}^{n-1} (3a^{d}_{234} - a^{d}_{244}) \end{array}.$$  Persevering a little, we work things out by hand for $n=1,2$ to find  $$\begin{array}{ccc}a_{225}^0 = 1 & a_{225}^1=4 & a_{225}^2=1 \\  a_{233}^0=0 & a_{233}^1=1 & a_{233}^2=4 \\ a_{234}^0=0 & a_{234}^1=5 & a_{234}^2=10 \\ a_{244}^0=1 & a_{244}^1=19 & a_{244}^2 = 25 \end{array}$$ $$\begin{array}{ccc} C_1^1=0 & C_1^2=-12 & C_1^3=0 \\ C_2^1=-6 & C_2^2 = -30 & C_2^3=0 \end{array}$$ so it follows that the invariants $a_{225}^3,a_{234}^3,a_{233}^3,a_{244}^3$ vanish because the RHS of the above system of equations is zero and the coefficient matrix is invertible.  Applying the same argument inductively (using $C_1^n = C_1^{n-1}+18a_{233}^{n-1}-6a_{234}^{n-1}$ and $C_2^n = C_2^{n-1}+18a_{234}^{n-1}-6a_{244}^{n-1}$) shows that all 3-point invariants for the class $(n,1)$ vanish when $n > 2$.

Next we compute the $3$-point invariants for $\beta = (n,2)$.  By the Dimension and Divisor Axioms it suffices to determine the $4$ numbers $$\begin{array}{cc} b^n_{255} := \langle T_2,T_5,T_5 \rangle^{n,2} & b^n_{335} := \langle T_3,T_3,T_5 \rangle^{n,2} \\ b^n_{345} := \langle T_3,T_4,T_5 \rangle^{n,2} & b^n_{445} := \langle T_4,T_4,T_5 \rangle^{n,2} \end{array}$$ because $\langle T_1,T_5,T_5 \rangle^{n,2} = (n/2)b^n_{255}$.  Apply WDVV with $\beta=(n,2)$ and $(i,j,k,l)=(1,4,3,3), \, (2,4,3,3), \, (2,3,4,4), \, (1,4,5,2)$.  In each of these cases, if $\beta_1 + \beta_2 = (n,2)$ and $\beta_1,\beta_2 \neq 0$ then each summand in the WDVV equation will vanish unless $\beta_1 = (d,1)$ and $\beta_2 = (n-d,1)$ for some $d \in \{0,\dots,n \}$.  This is because the $3$-point invariants for a homology class $(d,0)$ ($d>0$) vanish except $\langle T_1,T_1,T_1 \rangle^{d,0}$, but this never appears in the above WDVV equations because at most one of $i,j,k,l$ is $1$, so at most two insertions in any invariant are $T_1$.  Simplifying a little bit we find $$\begin{array}{rcl} b^n_{335} & = & \sum_{d=0}^n (2d-n)a_{233}^d a_{234}^{n-d} +d(n-d)(a_{234}^d a_{234}^{n-d}- a_{244}^d a_{233}^{n-d} ) \\
b^n_{345} & = & \sum_{d=0}^n (n-d)(3d-1)(a_{234}^d a_{234}^{n-d}- a_{244}^d a_{233}^{n-d} ) \\
b^n_{445} & = & \sum_{d=0}^n (3n-3d-1)(3d-1)(a_{234}^d a_{234}^{n-d}- a_{244}^d a_{233}^{n-d} )+(2d-n)a_{234}^d a_{244}^{n-d} \\
b^n_{255} & = & b_{445}^n + \sum_{d=0}^n d(3n-3d-1)a_{234}^d a_{225}^{n-d}-d(n-d) a_{244}^d a_{225}^{n-d}\end{array}$$ which implies that all these invariants vanish for $n>4$ (because of the vanishing of the $a^n$'s for $n>2$).  Working the rest out by hand, we find that they also vanish when $n=4$; the others are given below:  $$\begin{array}{cccc} b_{335}^0 = 0 & b_{335}^1 = 0 & b_{335}^2 = 6 & b_{335}^3 = 8 \\
b_{345}^0 = 0 & b_{345}^1 = 1 & b_{345}^2 = 20 & b_{345}^3 = 21 \\
b_{445}^0 = 0 & b_{445}^1 = 7 & b_{445}^2 = 64 & b_{445}^3 = 55 \\
b_{255}^0 = 0 & b_{255}^1 = 2 & b_{255}^2 = 8 & b_{255}^3 = 2 \end{array}$$  Previously, we gave a geometric reason for the vanishing of the $(0,2)$ invariants.

Finally we compute the $3$-point invariants $\langle T_5,T_5,T_5 \rangle^{n,3}$.  For dimension reasons, this is the only $3$-point invariant for $\beta=(n,3)$; when $\beta = (n,k)$ with $k>3$ then all $3$-point invariants vanish for dimension reasons.  To compute these invariants just apply WDVV with $\beta=(n,3)$ and $(i,j,k,l)=(2,3,5,5)$ to get: \begin{eqnarray*} \langle T_5,T_5,T_5 \rangle^{n,3} & = & \sum_{d=0}^n 3 \langle T_2,T_3,T_3 \rangle^{d,1} \langle T_1,T_5,T_5 \rangle^{n-d,2} - \langle T_2,T_3,T_4 \rangle^{d,1} \langle T_1,T_5,T_5 \rangle^{n-d,2} \\ & & - \langle T_2,T_3,T_3 \rangle^{d,1} \langle T_2,T_5,T_5 \rangle^{n-d,2}  - 3 \langle T_1,T_2,T_5 \rangle^{d,1} \langle T_3,T_3,T_5 \rangle^{n-d,2} \\ & & + \langle T_2,T_2,T_5 \rangle^{d,1} \langle T_3,T_3,T_5 \rangle^{n-d,2} +\langle T_1,T_2,T_5 \rangle^{d,1} \langle T_3,T_4,T_5 \rangle^{n-d,2}. \end{eqnarray*} Surely the RHS vanishes if $n>5$.  In fact, an explicit calculation shows that these invariants vanish when $n=0,1,5$ as well.  The rest of the cases can be easily computed from our previous results (see below for the values).
 
\bigskip

\centerline {\bf Quantum Cohomology of $Y$}

\bigskip

Here we assemble the 3-point invariants computed in the previous section to give a multiplication table and a presentation of the (small) quantum cohomology ring of $Y$.  This is a ring structure on $H^*(Y,\QQ) \otimes \QQ[q_2][[q_1]]$ with multiplication given by $$T_a \star T_b = T_aT_b+\sum_{(n,i) \neq (0,0)} \sum_{c=1}^5 \langle T_a,T_b,T_c \rangle^{n,i} T^c q_1^n q_2^i$$ (this is a finite sum for $(a,b) \neq (1,1)$).  The non-zero 3-point invariants of $Y$ for non-zero homology classes are given, up to reordering, by $\langle T_1,T_1,T_1 \rangle^{n,0}=6$ for $n>0$ and \begin{small} $$\begin{array}{cccc} \langle T_2,T_2,T_5 \rangle^{0,1}=1 & \langle T_2,T_4,T_4 \rangle^{0,1}=1 & \langle T_1,T_1,T_5 \rangle^{1,1}=4 &  \langle T_1,T_2,T_5 \rangle^{1,1}=4 \\  \langle T_1,T_3,T_3 \rangle^{1,1}=1 & \langle T_1,T_3,T_4 \rangle^{1,1}=5 &  \langle T_1,T_4,T_4 \rangle^{1,1}=19 &  \langle T_2,T_2,T_5 \rangle^{1,1}=4 \\ \langle T_2,T_3,T_3 \rangle^{1,1}=1 &  \langle T_2,T_3,T_4 \rangle^{1,1}=5 &
  \langle T_2,T_4,T_4 \rangle^{1,1}=19 &  \langle T_1,T_1,T_5 \rangle^{2,1}=4 \\  \langle T_1,T_2,T_5 \rangle^{2,1} = 2 &  \langle T_1,T_3,T_3 \rangle^{2,1}=8 &  \langle T_1,T_3,T_4 \rangle^{2,1}=20 &  \langle T_1,T_4,T_4 \rangle^{2,1}=50 \\  \langle T_2,T_2,T_5 \rangle^{2,1}=1 &  \langle T_2,T_3,T_3 \rangle^{2,1}=4 & \langle T_2,T_3,T_4 \rangle^{2,1}=10 & \langle T_2,T_4,T_4 \rangle^{2,1}=25 \end{array}$$ $$\begin{array}{cccc} 
\langle T_1,T_5,T_5 \rangle^{1,2}=1 & \langle T_1,T_5,T_5 \rangle^{2,2}=8  & \langle T_1,T_5,T_5 \rangle^{3,2}=3 \\
\langle T_2,T_5,T_5 \rangle^{1,2}=2 & \langle T_2,T_5,T_5 \rangle^{2,2}=8 & \langle T_2,T_5,T_5 \rangle^{3,2}= 2 \\
                                      & \langle T_3,T_3,T_5 \rangle^{2,2}=6 & \langle T_3,T_3,T_5 \rangle^{3,2}=8 \\
\langle T_3,T_4,T_5 \rangle^{1,2}=1 & \langle T_3,T_4,T_5 \rangle^{2,2}=20   & \langle T_3,T_4,T_5 \rangle^{3,2}=21 \\
\langle T_4,T_4,T_5 \rangle^{1,2}=7 & \langle T_4,T_4,T_5 \rangle^{2,2}=64  & \langle T_4,T_4,T_5 \rangle^{3,2}=55 \\ \langle T_5,T_5,T_5 \rangle^{2,3}=6 & \langle T_5,T_5,T_5 \rangle^{3,3} = 12 & \langle T_5,T_5,T_5 \rangle^{4,3}=6 \end{array}$$ \end{small} so the quantum multiplication table is as below. \begin{small} $$\begin{array}{rcl} T_1 \star T_1 & = & T_3+(-18T_3+6T_4)q_1(1-q_1)^{-1}+4q_1q_2+4q_1^2q_2 \\
T_1 \star T_2 & = & T_4+4q_1q_2+2q_1^2q_2 \\
T_1 \star T_3 & = & (2T_1+T_2)q_1q_2+(-4T_1+8T_2)q_1^2q_2 \\
T_1 \star T_4 & = & T_5 + (4T_1+5T_2)q_1q_2+(-10T_1+20T_2)q_1^2q_2 \\
T_1 \star T_5 & = & (-8T_3+4T_4)q_1q_2+(-10T_3+4T_4)q_1^2q_2 + q_1q_2^2+8q_1^2q_2^2+3q_1^3q_2^2 \\
T_2 \star T_2 & = & -6T_3+3T_4+q_2+4q_1q_2+q_1^2q_2 \\
T_2 \star T_3 & = & T_5+(2T_1+T_2)q_1q_2+(-2T_1+4T_2)q_1^2q_2 \\
T_2 \star T_4 & = & 3T_5+T_1q_2+(4T_1+5T_2)q_1q_2+(-5T_1+10T_2)q_1^2q_2 \\
T_2 \star T_5 & = & T_3q_2+(-8T_3+4T_4)q_1q_2+(-5T_3+2T_4)q_1^2q_2 + 2q_1q_2^2+8q_1^2q_2^2+2q_1^3q_2^2 \\
T_3 \star T_3 & = & (-2T_3+T_4)q_1q_2 + (-20T_3+8T_4)q_1^2q_2 + 6q_1^2q_2^2 + 8 q_1^3q_2^2 \\
T_3 \star T_4 & = & (-10T_3+5T_4)q_1q_2+(-50T_3+20T_4)q_1^2q_2 +q_1q_2^2+20q_1^2q_2^2 + 21q_1^3q_2^2 \\
T_3 \star T_5 & = & T_1q_1q_2^2+(2T_1+6T_2)q_1^2q_2^2+(-3T_1+8T_2)q_1^3q_2^2 \\
T_4 \star T_4 & = & T_3q_2+(-38T_3+19T_4)q_1q_2+(-50T_3+50T_4)q_1^2q_2+7q_1q_2^2+64q_1^2q_2^2+55q_1^3q_2^2 \\
T_4 \star T_5 & = & (4T_1+T_2)q_1q_2^2+(4T_1+20T_2)q_1^2q_2^2+(-8T_1+21T_2)q_1^3q_2^2 \\
T_5 \star T_5 & = & (-T_3+T_4)q_1q_2^2+(-16T_3+8T_4)q_1^2q_2^2 + (-7T_3+3T_4)q_1^3q_2^2 \\ & & + 6q_1^2q_2^3+12q_1^3q_2^3+6q_1^4q_2^3 \end{array}$$  \end{small} Deforming the relations in $H^*(Y,\QQ)$, we can give a presentation: $$QH^*(Y) \cong \QQ[T_1,T_2,q_2][[q_1]]/ \langle R_1, R_2 \rangle$$ where \begin{small} $$\begin{array}{l} R_1 = T_2^2-3T_1T_2+6T_1^2-q_2-16q_1q_2-19q_1^2q_2+18q_1(1-q_1)^{-1}(T_1T_2-T_2^2+q_2-q_1^2q_2) \\ R_2 = T_1^3+(-6T_1-T_2)q_1q_2-8T_2q_1^2q_2-6q_1(1-q_1)^{-1}[T_1^2T_2+(-10T_1-3T_2)q_1q_2 \\ \quad \quad \quad + (10T_1-24T_2)q_1^2q_2]. \end{array}$$ \end{small}

\end{section}

\begin{section}{The Change of Variables} \label{section:changevars}

Since $QH_{orb}^*([F/\ZZ_2])$ has no degree $0$ quantum parameter, we begin by setting the quantum parameter $q_1$ to $-1$, in which case $QH^*(Y)$ is determined by the part of the multiplication table below $$\begin{array}{lcl} T_1 \star T_1 = 10T_3-3T_4 & \quad & T_2 \star T_2 = 3T_4-6T_3-2q_2 \\ T_1 \star T_2 = T_4-2q_2 & \quad & T_2 \star T_3 = T_5+(-4T_1+3T_2)q_2 \\ T_1 \star T_3 = (-6T_1+7T_2)q_2 & \quad & T_2 \star T_4 = 3T_5+(-8T_1+5T_2)q_2 \\ T_1 \star T_4 = T_5 + (-14T_1+15T_2)q_2 \end{array}$$ and the presentation simplifies similarly: $$QH^*(Y)|_{q_1=-1} \cong \QQ[T_1,T_2,q_2]/ \langle \tilde{R_1}, \tilde{R_2} \rangle$$ where $$\begin{array}{l} \tilde{R_1} = 5T_2^2-6T_1T_2+3T_1^2-2q_2 \\ \tilde{R_2} = T_1^3+3T_1^2T_2+(66T_1-70T_2)q_2. \end{array}$$ 

Now it is straightforward to check that $QH_{orb}^*( \X ) \otimes_\QQ \QQ(i)$ and $QH^*(Y)|_{q_1=-1} \otimes_\QQ \QQ(i)$ are isomorphic by the maps below.  $$\begin{array}{ccc} S_1 & \mapsto & T_2 \\ S_2 & \mapsto & i(T_2-T_1) \\ q & \mapsto & -q_2 \end{array} \quad \quad \begin{array}{ccc} T_1 & \mapsto & S_1+iS_2 \\ T_2 & \mapsto & S_1 \\ q_2 & \mapsto & -q \end{array}$$  This is the same change of variables used in \cite{W}.  Clearly these are inverse maps, so all we need to do is show that they are well defined (i.e. that they kill the relations).  The most difficult such computation is checking that $S_2^3-6qS_2 \mapsto 0$.  To do this, use the relations $T_1 \tilde{R_1}, T_2 \tilde{R_1}$, and $\tilde{R_2}$ (or the multiplication table) to express the degree three monomials in terms of $T_1^2T_2$: $$\begin{array}{lll} T_1^3 & = & -3T_1^2T_2+(-66T_1+70T_2)q_2 \\ T_1T_2^2 & = & 3T_1^2T_2+(40T_1-42T_2)q_2 \\ T_2^3 & = &  3T_1^2T_2 + (48T_1-50T_2)q_2 \end{array}$$ so that $$\begin{array}{lll}-i(S_2^3-6qS_2) & \mapsto & T_1^3-3T_1^2T_2+3T_1T_2^2-T_2^3+6q_2(T_2-T_1) \\ & = & -3T_1^2T_2+(-66T_1+70T_2)q_2 \\ &  & -3T_1^2T_2 \\ &  & +9T_1^2T_2+(120T_1-126T_2)q_2  \\ & & -3T_1^2T_2+(-48T_1+50T_2)q_2 + (-6T_1+6T_2)q_2 \\ & = & 0.\end{array}$$ The other checks are similar, though less difficult.

\end{section}

\begin{section}{Higher Point Invariants of $\X$} \label{section:higherpoint}

In this section we show that the genus zero Gromov-Witten invariants of $\X$ are determined by the $3$-point invariants that we already computed (in fact we will show that only the $2$-point invariants are needed) together with the invariants $$\langle S_2,\dots,S_2 \rangle^0,$$ which we will compute below.  In what follows, by a \emph{divisor}, we will mean a ``non-twisted" divisor: a cohomology class $\delta \in H^2( \X ) \subseteq H^2_{orb}(\X)$.  Recall that the usual WDVV equation \ref{wdvv} holds for orbifolds without modification \cite{AGV2}, but the Point Mapping and Divisor Axioms do not (hence \ref{lhs} and \ref{rhs} are not quite right for orbifolds, though the extra terms that show up in \ref{lhs} and \ref{rhs} occur when $0 < |V| < n$ and hence only involve invariants of lower point number).  This is because both of these axioms are proved (for smooth varieties) using the forgetful stabilization map $$\ov{M}_{0;n+1}(X,\beta) \to \ov{M}_{0;n}(X,\beta)$$ (which exists when $n \geq 3$ or $\beta \neq 0$).  However, when $\X$ is an orbifold (say $\X = [X / \ZZ_2]$) then, although there is a forgetful map $$\ov{M}_{0;r,u+1}(\X,\beta) \to \ov{M}_{0;r,u}(\X,\beta),$$ (at least when $r+u \geq 3$ or $\beta \neq 0$) there is \emph{no} forgetful map $$\ov{M}_{0;r+1,u}(\X,\beta) \to \ov{M}_{0;r,u}(\X,\beta).$$  

Now we recall a simple reconstruction theorem for genus zero Gromov-Witten invariants, which is particularly relevant for orbifolds.  This is essentially the natural statement of Kontsevich's reconstruction theorem when the cohomology ring is not generated by divisors; it appears in \cite{W} in essentially the same form as below (we include the proof here as well since it is so simple).  Let $H^*_{orb}(\X)_{div}$ denote the subring of $H^*_{orb}(\X)$ generated by divisors (in the sense above).  

\begin{lem} \label{lem:reconstruction} Suppose that $\theta_1,\dots,\theta_N \in H^*_{orb}(\X)$ generate $H^*_{orb}(\X)$ as a module over $H^*_{orb}(\X)_{div}$.  Then the genus zero Gromov-Witten invariants of $\X$ are determined by linearity, the WDVV equations, the Divisor Axiom, the $2$-point invariants, the $3$-point degree $0$ invariants, and the invariants of the form $$\langle \eta \theta_{i_1},\theta_{i_2},\dots,\theta_{i_k} \rangle^\beta,$$ where $\eta \in H^*_{orb}(\X)_{div}$.
\end{lem}

\begin{proof} We say $\langle \alpha_1, \dots, \alpha_n \rangle^{\beta_1}$ is \emph{lower} than $\langle \alpha^{'}_1, \dots, \alpha^{'}_m \rangle^{\beta_2}$ if $n \leq m$, $\beta_1 \leq \beta_2$, and one of these inequalities is strict.  By the hypothesis on the $\theta_i$ and linearity it suffices to show that an invariant of the form $$\langle \delta_{1,i_1} \delta_{1,i_2} \cdots \delta_{1,i_{m_1}} \theta_{j_1} , \delta_{2,i_1} \delta_{2,i_2} \cdots \delta_{2,i_{m_2}} \theta_{j_2} , \dots, \delta_{n,i_1} \delta_{n,i_2} \cdots \delta_{n,i_{m_n}} \theta_{j_n} \rangle^\beta,$$ where the $\delta_{i,j}$ are various divisors depending on the pair $(i,j)$ (some $m_j$ may be zero), is determined by the data mentioned above.  We may assume by induction that all lower invariants are determined by such data and that $n \geq 3$.  Notice that if $\delta$ is a divisor and $$\langle \alpha_1,\dots,\alpha_{i-1}, \delta \alpha_i, \alpha_{i+1} , \dots , \alpha_n \rangle^\beta$$ is an invariant with $n \geq 3$, then applying WDVV (\ref{wdvv}) with $\phi_1=\alpha_i$, $\phi_2 = \delta$, $\phi_3=\alpha_1$, and $\phi_4,\gamma_1,\dots,\gamma_{n-3}$ equal to $\alpha_2,\dots,\alpha_{i-1},\alpha_{i+1},\dots,\alpha_n$ expresses the above invariant in terms of lower invariants, the invariant $$\langle \delta \alpha_1, \alpha_2, \dots, \alpha_n \rangle^\beta ,$$ and two other invariants which are immediately reduced to lower invariants by the Divisor Axiom.  Applying this repeatedly we can move all the $\delta_{i,j}$ to the first insertion to establish the lemma. \end{proof}

When $\X=[F / \ZZ_2]$, $\theta_1=S_2$ satisfies the hypotheses on the $\theta_i$ above.  In $H^*_{orb}(\X)$, $S_1S_2=4S_4$, $S_1^n S_2 = 0$ for $n > 1$, and the invariants $$\langle S_4,S_2,S_2,\dots,S_2 \rangle^d$$ all vanish for dimension reasons, so the genus $0$ Gromov-Witten theory is determined by previously computed invariants, together with the invariants of the form $$\langle  \underbrace{ S_2,\dots,S_2 }_{2g+2} \rangle^0.$$  The corresponding moduli space $\ov{M}_{0;2g+2,0}(\X,0)$ is just $C \times \ov{M}_{0;2g+2,0}( \B \ZZ_2)$ and the virtual fundamental class is dual to the Euler class of $$\pi_1^* N_{C / F} \otimes \pi_2^* \EE^\lor,$$ where $\EE$ is the Hodge bundle on $\ov{M}_{0;2g+2,0}( \B \ZZ_2)$.  Let $\lambda_i := c_i( \EE )$.  Since $p_1|_C=p_2|_C=2$, and $c_1(TF)=2p_1+2p_2$, the SES $$0 \to TC \cong T \PP^1 \to TF|_C \to N_{C / F} \to 0$$ shows that $c_1( N_{C / F} ) = 6$.  Thus we compute \begin{eqnarray*} \langle \underbrace{S_2,\dots,S_2}_{2g+2} \rangle^0 & = & \int_{ \ov{M}_{0;2g+2,0}(\X,0) } e( {\rm Obstruction \; Bundle} ) \\ & = & \int_{C \times \ov{M}_{0;2g+2,0}( \B \ZZ_2)} c_{2g}( \pi_1^* N_{C / F} \otimes \pi_2^* \EE^\lor ) \\ & = & -6 \int_{\ov{M}_{0;2g+2,0}( \B \ZZ_2)} \lambda_g \lambda_{g-1} \end{eqnarray*}  These integrals were computed in \cite{FP} and are well-known in the subject.

To complete the proof of the Crepant Resolution Conjecture, we determine the change of variables to use by setting $q=-q_2=0$ in Section~\ref{section:changevars} and determine the matrix giving the change of variables $\{ S_i \} \to \{ T_i \}$.  Thinking for a moment about how the potential function encodes the product in these isomorphic rings, we see that we should use the transpose of this matrix as the change of variables $\{ t_i \} \to \{ s_i \}$ for the potential functions.  Indeed, we have shown that the potential functions of $Y$ and $\X$ are given as below (neglecting terms of degree $\geq 4$ in the $s_i$ and $t_i$, as well as ``unstable terms" of degree $<3$). \begin{small} \begin{eqnarray*} F^Y(t_0,\dots,t_5,q_1,q_2) & = & \frac{t_0^2t_5}{2}+t_0t_1t_4+t_0t_2t_3+t_0t_2t_4+\frac{t_1^2t_2}{2}+\frac{3t_1t_2^2}{2}+\frac{t_2^3}{2}+\frac{t_1^3q_1}{1-q_1} \\ & & +\frac{t_2^2t_5q_2}{2}+\frac{t_2t_4^2q_2}{2}+2t_1^2t_5q_1q_2+4t_1t_2t_5q_1q_2 + \frac{t_1t_3^2q_1q_2}{2}\\ & & +5t_1t_3t_4q_1q_2+\frac{19t_1t_4^2q_1q_2}{2}+2t_2^2t_5q_1q_2+\frac{t_2t_3^2q_1q_2}{2}+5t_2t_3t_4q_1q_2 \\ & & + \frac{19 t_2t_4^2q_1q_2}{2} + 2 t_1^2t_5q_1^2q_2+2t_1t_2t_5q_1^2q_2+4t_1t_3^2q_1^2q_2 \\ & & + 20t_1t_3t_4q_1^2q_2+25t_1t_4^2q_1^2q_2+\frac{t_2^2t_5q_1^2q_2}{2}+2t_2t_3^2q_1^2q_2+10t_2t_3t_4q_1^2q_2 \\ & & + \frac{25t_2t_4^2q_1^2q_2}{2}+\frac{t_1t_5^2q_1q_2^2}{2} + 4t_1t_5^2q_1^2q_2^2+\frac{3t_1t_5^2q_1^3q_2^2}{2} \\ & & + t_2t_5^2q_1q_2^2+4t_2t_5^2q_1^2q_2^2+t_2t_5^2q_1^3q_2^2+3t_3^2t_5q_1^2q_2^2+4t_3^2t_5q_1^3q_2^2 \\ & & t_3t_4t_5q_1q_2^2+20t_3t_4t_5q_1^2q_2^2+21t_3t_4t_5q_1^3q_2^2+\frac{7t_4^2t_5q_1q_2^2}{2} \\ & & + 32t_4^2t_5q_1^2q_2^2+\frac{55t_4^2t_5q_1^3q_2^2}{2}+t_5^3q_1^2q_2^3+2t_5^3q_1^3q_2^3+t_5^3q_1^4q_2^3+O(4) \end{eqnarray*} \begin{eqnarray*} F^\X(s_0,\dots,s_5,q) & = & \frac{3s_0^2s_5}{2}+3s_0s_1s_3+\frac{s_0s_2s_4}{2}+\frac{s_1^3}{2}+s_1s_2^2+\frac{9s_1s_3^2q}{2}+3s_1^2s_5q \\ & &+27s_3^2s_5q^2+18s_1s_5^2q^2+\frac{s_1s_4^2q}{2}+3s_2s_3s_4q+3s_2^2s_5q+\frac{3s_4^2s_5q^2}{2} \\ & & + O(4) \end{eqnarray*} \end{small}

It is easy (with a computer, say) to check that $$F^Y(s_0,-is_2,s_1+is_2,-6s_3-\frac{3 i }{2}s_4,3s_3+is_4,3s_5,-1,-q) = F^{\X}(s_0,\dots,s_5,q).$$  Since the change of variables is linear and respects the metrics, it preserves WDVV equations, so the entire potentials will agree (because of Lemma~\ref{lem:reconstruction}) if the coefficients of $s_2^n$ (with no powers of $q$) agree.  In particular, it will be sufficient to show that the full potential functions agree (up to unstable terms) under the above change of variables (analytically continuing to $q_1=-1$) when we set $q=s_0=s_3=s_4=s_5=0$.  Notice that, because of the Divisor Axiom and the fact that $T_2$ evaluates $0$ on a homology class of the form $(l,0)$, the dependence on $t_2$ in the potential function of $Y$ is purely classical when we set $q_2=0$.  Also notice that the coefficient of $t_1^n$ (when $q_2=0$) is determined by the invariant $\langle T_1,\dots,T_1 \rangle^{d,0}$, which reduces by the Divisor Axiom to our previous computation of $\langle \rangle^{d,0}$.  Putting these observations together, we have \begin{eqnarray*} F^Y(0,t_1,t_2,0,0,0,q_1,0) & = & \frac{t_1^2t_2}{2}+\frac{3t_1t_2^2}{2}+\frac{t_2^3}{2}+6 \sum_{d \geq 1} \frac{1}{d^3} e^{dt_1}q_1^d \end{eqnarray*} and, using the computation above, $$F^{\X}(0,s_1,s_2,0,0,0,0)  =  \frac{s_1^3}{2}+s_1s_2^2-6 \sum_{g \geq 1} \frac{1}{(2g+2)!} \left ( \int_{ \ov{M}_{0;2g+2,0}( \B \ZZ_2) } \lambda_g \lambda_{g-1} \right ) s_2^{2g+2},$$  so $$F^Y(0,-is_2,s_1+is_2,0,0,0,0)  =  \frac{s_1^3}{2}+s_1s_2^2 + 6 \left ( \frac{i s_2^3}{12} + \sum_{d \geq 1} \frac{1}{d^3} e^{-ids_2} q_1^d \right ).$$ Clearly all the terms involving $s_1$ match up, and the fact that the third derivatives with respect to $s_2$ agree after analytically continuing to $q_1=-1$ is the same computation made for $[\CC^2 / \ZZ_2]$ in \cite{BG}, except with a factor of $6$ floating around\footnote{Also $i$ is replaced by $-i$, but that makes no difference since the result is real anyway.} (both of these third partials are equal to $-3 \tan (s_2 / 2)$).

\end{section}

\end{document}